\documentclass[12pt]{amsart}
\usepackage{amsthm}
\usepackage{amsmath, mathrsfs}
\usepackage{color}
\usepackage{amssymb}
\usepackage{hyperref}
\usepackage{enumerate}
\usepackage{latexsym,array}
\usepackage{amsfonts}
\usepackage{shadow}
\newcommand{\ran}{{\rm ran}\,}

\newtheorem{Pa}{Paper}[section]
\newtheorem{Tm}[Pa]{{\bf Theorem}}
\newtheorem{La}[Pa]{{\bf Lemma}}
\newtheorem{Dn}[Pa]{{\bf Definition}}

\newtheorem{Rk}[Pa]{{\bf Remark}}
\newtheorem{Pn}[Pa]{{\bf Proposition}}

\def\e{\epsilon_N}

\def\hh{\mathbb{H}}
\author[D. Alpay]{Daniel Alpay}
\address{(DA) Department of Mathematics\\
BenGurion University of the Negev\\
Beer-Sheva 84105 Israel} \email{dany@math.bgu.ac.il}

\author[F. Colombo]{Fabrizio Colombo}
\address{(FC) Politecnico di
Milano\\Dipartimento di Matematica\\Via E. Bonardi, 9\\20133
Milano, Italy}
\email{fabrizio.colombo@polimi.it}
\author[I. Sabadini]{Irene Sabadini}
\address{(IS) Politecnico di
Milano\\Dipartimento di Matematica\\Via E. Bonardi, 9\\20133
Milano, Italy}
\email{irene.sabadini@polimi.it}
\title[Krein-Langer factorization]
{Krein-Langer factorization and related topics in the slice
hyperholomorphic setting} \oddsidemargin 0.2in \evensidemargin 0.2in
\def\s{\sigma}

 \keywords{Schur functions,
realization, reproducing kernels, slice hyperholomorphic
functions, $S$-resolvent operators.}


\thanks{D. Alpay thanks the Earl Katz family for endowing the chair
which supported his research, and the Binational Science
Foundation Grant number 2010117.  F. Colombo and I. Sabadini acknowledge the Center for Advanced Studies of the Mathematical Department of the Ben-Gurion University of the Negev for the support and the kind hospitality during the period in which this paper has been written.}
\begin{document}
\maketitle \tableofcontents
\parindent 0cm
\begin{abstract}
We study various aspects of Schur analysis in the slice
hyperholomorphic setting.  We present two sets of results: first, we
give new results on the functional calculus for slice
hyperholomorphic functions. In particular, we introduce and study some properties of the Riesz projectors.   Then we prove a Beurling-Lax type
theorem, the so-called structure theorem. A crucial fact which allows to prove our results, is the fact that
the right spectrum of a quaternionic linear operator and the point S-spectrum coincide.
Finally, we study the Krein-Langer factorization
for slice hyperholomorphic generalized Schur functions. Both the Beurling-Lax type theorem and the Krein-Langer factorization are far reaching results which have not been proved in the quaternionic setting using notions of hyperholomorphy other than slice hyperholomorphy.
\end{abstract}
\section{Introduction}
\setcounter{equation}{0}
In operator theory and in the theory of linear systems, the study
of functions analytic and contractive in the open unit disk
(Schur functions) is called Schur analysis. It includes, in
particular, interpolation problems, operator models, and has been
extended to various settings. See for instance
\cite{MR2002b:47144,c-schur,hspnw,goh1} for some related books.
In \cite{acs1,acs2} we began a study of Schur analysis in the
setting of slice hyperholomorphic functions.
Following \cite{acs1,acs2}, let us recall that a
generalized Schur function is a $\mathbb H^{N\times M}$ valued
function $S$ slice-hyperholomorphic in a neighborhood $V$ of the
origin and for which the kernel
\begin{equation}
\label{skl}
K_S(p,q)= \sum_{n=0}^\infty p^n(I_N-S(p)
S(q)^*)\overline{q}^n
\end{equation}
has a finite number of negative squares in $V$, or more generally
such that the kernel
\begin{equation}
\label{kernpot}
\sum_{n=0}^\infty p^n(\s_2-S(p)
\s_1S(q)^*)\overline{q}^n
\end{equation}
has a finite number of negative squares in $V$, where
$\s_1\in\mathbb H^{M\times M}$ and $\s_2\in\mathbb H^{N\times N}$
are signature matrix (i.e. both self-adjoint and invertible).
Since this work is aimed at different audiences, it is worth
mentioning that the classical counterparts of the kernels
\eqref{kernpot} originate with the theory of characteristic
operator functions. In the indefinite case, such kernels have
been studied by Krein and Langer; see for instance
\cite{kl1,MR47:7504}. When $\sigma_1=\sigma_2$ and when the
kernel is positive definite, Potapov gave in the fundamental paper
\cite{pootapov} the multiplicative structure of the corresponding
functions $S$.\\

In \cite{acs1} we studied the realization of such $S$ in the case
of positive definite kernels. In \cite{acs2} we studied an
interpolation problem, and began the study of the indefinite
metric case, where the underlying spaces are Pontryagin spaces
rather than Hilbert spaces. In this work we prove a Beurling-Lax
type theorem in this setting and study the Krein-Langer
factorization for slice hyperholomorphic generalized Schur
functions. Slice hyperholomorphic functions turned out to be a
natural tool to generalize Schur analysis to the quaternionic
setting. Some references for this theory of functions, with no
claim of completeness, are \cite{MR2353257, MR2555912, MR2742644},
the book \cite{MR2752913} and the forthcoming \cite{GSS}.\\

The analogue of the
resolvent operator in classical analysis is now the $S$-resolvent
 operator, and according to this resolvent, the spectrum has to
be replaced by the $S$-spectrum. The relation between the
$S$-spectrum and the right spectrum of a right linear
quaternionic operator  is important for the present paper.
Indeed, in the literature there are several results on the right
spectrum which is widely studied, especially for its application
in mathematical physics, see e.g. \cite{MR1333599}. However, it
is well known that the right spectrum is not associated to a right
linear quaternionic operator; the eigenvectors associated to a
given eigenvalue do not even form a linear space. The
$S$-spectrum  arises in a completely different setting, it is
associated to a right linear operator and, quite surprisingly, the point $S$-spectrum 
coincides with the right spectrum. This fact and the fact that
any right eigenvector is also an S-eigenvector, see Proposition
\ref{eigenvector}, allow to use for the point $S$-spectrum various
result which hold for the right spectrum, see Sections 6 to 9.\\

The  $S$-resolvent operator  allows the
definition of the quaternionic analogue of the operator
$(I-zA)^{-1}$ that appears in the realization function
$s(z)=D+zC(I-zA)^{-1}B$. It turns out that when $A$ is a
quaternionic matrix and $p$ is a quaternion then $(I-pA)^{-1}$
has to be replaced by $( I  -\bar p A)(|p|^2A^2-2{\rm Re}(p) A+
I  )^{-1}$ which is equal to $p^{-1}S^{-1}_R(p^{-1},A)$ where
$S^{-1}_R(p^{-1},A)$ is the right $S$-resolvent operator
associated to the quaternionic matrix $A$.
Moreover, the $S$-resolvent operator allows to introduce and study the Riesz projectors and
 the invariant subspaces under a quaternionic operator.
\\

The S-resolvent operator is also a fundamental tool to define
the quaternionic functional calculus, and we refer the reader to
\cite{MR2735309, MR2496568, MR2661152, MR2803786} for further discussions.
Schur multipliers in the quaternionic setting have been studied also
in \cite{MR2124899,MR2240272,asv-cras},
in a different setting, using the
Cauchy-Kovalesvkaya product and series of Fueter polynomials.
Since Schur analysis plays an important role in linear systems, we mention
that papers \cite{MR733953,MR2205693,MR2131920} treat
various aspects of a theory of linear systems in the quaternionic
setting.
We finally remark that it is possible to define slice hyperholomorphic functions with values
in a Clifford algebra, \cite{MR2520116, MR2684426, MR2673423}, which admit a functional
calculus for $n$-tuples of
operators, see \cite{MR2402108, MR2720712, SCFUNCTIONAL,MR2752913}.
\\

The paper consists of eight sections besides the introduction,
and its outline is as follows: Sections \ref{2}-\ref{4} are related
to results on slice hyperholomorphic functions and the related
functional calculus. Sections \ref{6}-\ref{9} are related to Schur
analysis. More precisely: Sections \ref{2} and \ref{3} are of a
survey nature on slice hyperholomorphic functions and the quaternionic
functional calculus, respectively. Section \ref{4} contains new
results on the analogue of Riesz projector for the quaternionic
functional calculus. Moreover, it contains a discussion on the
right spectrum, which has been widely studied in the literature
both in linear algebra \cite{MR97h:15020} and in mathematical physics
\cite{MR1333599}, and which, as we have already pointed out,
coincides with the point $S$-spectrum. These results will be used
 in the second part of the paper. A characterization of the number of
negative squares of a slice hyperholomorphic kernel in terms of
its power series coefficients is given in Section \ref{6}. In
Section \ref{7} we present some results on linear operators in
quaternionic Pontryagin spaces. We show in particular that a
contraction with no S-spectrum on the unit sphere has a unique
maximal negative invariant subspace. In Section \ref{8} we prove
a version of the Beurling-Lax theorem, the so-called structure
theorem, in the present setting. In Section \ref{5} we discuss
the counterparts of matrix-valued unitary rational functions. The last section considers a
far reaching result in the quaternionic framework, namely the
Krein-Langer factorization theorem for generalized Schur
functions. It is interesting to note that the result is based on
Blaschke products whose zeros and poles have a peculiar behaviour
when taking the slice hyperholomorphic reciprocal.

\section{Slice hyperholomorphic functions}
\setcounter{equation}{0}
\label{2}
We begin this section by recalling the notations and some basic facts on the theory of slice
hyperholomorphic functions that we will use in the sequel. We send the reader to the papers
\cite{MR2555912, MR2742644, MR2353257} and the book \cite{MR2752913}
for more details. Let
$\hh$ be the real associative algebra of quaternions with respect
to the basis $\{1, i,j,k \}$ whose elements satisfy the relations
$$
i^2=j^2=k^2=-1,\
 ij =-ji =k,\
jk =-kj =i ,
 \  ki =-ik =j .
$$
We will denote a quaternion $p$ as $p=x_0+ix_1+jx_2+kx_3$,
$x_\ell\in \mathbb{R}$, $\ell=0,1,2,3$, its conjugate as
$\bar p=x_0-ix_1-jx_2-kx_3$, its norm $|p|^2=p\overline p$.
The real part of a quaternion will be denoted with the symbols ${\rm Re}(p)$ or $x_0$, while ${\rm Im}(p)$ denotes the imaginary part of $p$.
\noindent
Let $\mathbb{S}$ be the 2-sphere of purely imaginary unit quaternions, i.e.
\begin{equation}
\label{sphere}
\mathbb{S}=\{ p=ix_1+jx_2+kx_3\ |\
x_1^2+x_2^2+x_3^2=1\}.
\end{equation}
To each nonreal quaternion $p$ it is possible to uniquely associate the
element $I_p\in\mathbb{S}$ defined by
$$
I_p=
\displaystyle\frac{{\rm Im}(p)} {|{\rm Im}(p)|}.
$$
The complex plane $\mathbb{C}_{I_p}=\mathbb{R}+I_p\mathbb{R}=\{x+I_qy \ | \  x,y\in \mathbb{R}\}$ is determined by the imaginary unit $I_p$, and $\mathbb{C}_{I_p}$ obviously contains $p$.
\begin{Dn}
Given $p\in\hh$, $p=p_0+I_pp_1$ we denote by $[p]$ the set of all
elements of the form $p_0+Jp_1$ when $J$ varies in $\mathbb{S}$.
\end{Dn}
\noindent
\begin{Rk}{\rm
The set $[p]$ is a $2$-sphere which is reduced to the point $p$
when $p\in\mathbb{R}$.}
\end{Rk}
\noindent
We now recall the definition of slice hyperholomorphic functions.
\begin{Dn}[Slice hyperholomorphic functions]
Let $U\subseteq\hh$ be an open set and let
$f:\ U\to\hh$ be a real differentiable function. Let
$I\in\mathbb{S}$ and let $f_I$ be the restriction of $f$ to the
complex plane $\mathbb{C}_I := \mathbb{R}+I\mathbb{R}$ passing through $1$
and $I$ and denote by $x+Iy$ an element on $\mathbb{C}_I$.
\begin{itemize}
\item[(1)]
 We say that $f$ is a left slice hyperholomorphic function
(or hyperholomorphic for short)  if, for every
$I\in\mathbb{S}$, we have
$$
\frac{1}{2}\left(\frac{\partial }{\partial x}+I\frac{\partial
}{\partial y}\right)f_I(x+Iy)=0.
$$
\item[(2)]
We say that $f$ is right slice hyperholomorphic function (or right
hyperholomorphic for short) if,
for every
$I\in\mathbb{S}$, we have
$$
\frac{1}{2}\left(\frac{\partial }{\partial x}f_I(x+Iy)+\frac{\partial
}{\partial y}f_I(x+Iy)I\right)=0.
$$
\item[(3)]
In the sequel we will denote by $\mathcal{R}^L(U)$ (resp.
$\mathcal{R}^R(U)$)  the right (resp. left)  $\mathbb{H}$-vector space of
left (resp. right) hyperholomorphic functions on the open set $U$.
When we do not distinguish between $\mathcal{R}^L(U)$
and $\mathcal{R}^R(U)$ we will use the symbol
$\mathcal{R}(U)$.
\end{itemize}
\end{Dn}
The natural open sets on which slice hyperholomorphic functions are defined are axially symmetric, i.e. open sets that contain the 2-sphere
$[p]$ whenever they contain $p$, which are also s-domains, i.e. they are domains which remain connected when intersected with any complex plane $\mathbb{C}_I$.
\\
Given two left slice
hyperholomorphic functions $f$, $g$, it is possible to introduce
a binary operation called the $\star$-product, such that $f\star g$ is a slice
hyperholomorphic function.
Let $f,g:\ \Omega \subseteq\mathbb{H}$ be slice hyperholomorphic functions such that
their restrictions to the complex plane $\mathbb{C}_I$ can be written as
$f_I(z)=F(z)+G(z)J$,
$g_I(z)=H(z)+L(z)J$ where $J\in\mathbb{S}$, $J\perp I$. The functions $F$,
$G$, $H$, $L$ are holomorphic functions of the variable $z\in
\Omega \cap \mathbb{C}_I$ and they exist by the splitting lemma, see \cite[p. 117]{MR2752913}.
We can now give the following:
\begin{Dn}
Let $f,g$ slice hyperholomorphic functions defined on an axially symmetric open set $\Omega\subseteq\mathbb{H}$.
The $\star$-product of  $f$ and $g$ is defined as the unique
left slice hyperholomorphic function on $\Omega$ whose restriction to the
complex plane $\mathbb{C}_I$ is given by
\begin{equation}\label{starproduct}
(F(z)+G(z)J)\star(H(z)+L(z)J):=
(F(z)H(z)-G(z)\overline{L(\bar z)})+(G(z)\overline{H(\bar z)}+F(z)L(z))J.
\end{equation}
\end{Dn}
When $f$ are expressed by power series, i.e. $f(p)=\sum_{n=0}^\infty p^n a_n$, $g(p)=\sum_{n=0}^\infty p^n b_n$, then $(f\star g)(p)=\sum_{n=0}^\infty p^n c_n$ where
$c_n=\sum_{r=0}^na_rb_{n-r}$ is obtained by convolution on the coefficients. This product extends the product of quaternionic polynomials with right coefficients, see \cite{lam}, to series.
Analogously, one can introduce a $\star$-product for right slice
hyperholomorphic functions. For more details we refer the reader to \cite{MR2752913}. When considering in a same
formula both the products, or when confusion may arise, we will write $\star_l$ or
$\star_r$ according to the fact that we are using the left or the
right slice hyperholomorphic product. When there is no subscript, we will
mean that we are considering the left $\star$-product.\\

 Given a slice hyperholomorphic function $f$, we
can define its slice hyperholomorphic reciprocal $f^{-\star}$, see
\cite{MR2555912, MR2752913}. In this paper it will be sufficient to
know the following definition
\begin{Dn}\label{reciprocal}
Given $f(p)=\sum_{n=0}^\infty p^n a_n$, let
us set
$$
f^c(p)=\sum_{n=0}^\infty p^n \bar a_n,\qquad  f^s(p)=(f^c\star f)(p
)=\sum_{n=0}^\infty p^nc_n,\quad
c_n=\sum_{r=0}^n a_r\bar a_{n-r},
$$
where the series converge.
The left slice hyperholomorphic reciprocal of $f$
is then defined as
$$
f^{-\star}:=(f^s)^{-1}f^c.
$$
\end{Dn}

\section{Formulations of the quaternionic functional calculus}
\label{3}
Here we briefly recall the possible formulations of the quaternionic functional calculus
that we will use in the sequel.
Let $V$ be a two sided quaternionic Banach space, and let
$\mathcal{B}(V)$ be the two sided vector space of all right
linear bounded operators on $V$.

\begin{Dn}[The $S$-spectrum and the $S$-resolvent sets of
quaternionic operators]\label{defspscandres}
Let $T\in\mathcal{B}(V)$.
We define the $S$-spectrum $\sigma_S(T)$ of $T$  as:
$$
\sigma_S(T)=\{ s\in \mathbb{H}\ \ :\ \ T^2-2 {\rm Re}\,(s) T+|s|^2\mathcal{I}\ \ \
{\it is\ not\  invertible}\}.
$$
The $S$-resolvent set $\rho_S(T)$ is defined by
$$
\rho_S(T)=\mathbb{H}\setminus\sigma_S(T).
$$
\end{Dn}

The notion of $S$-spectrum of a linear quaternionic operator $T$ is suggested
by the definition of $S$-resolvent operator that is the analogue of the Riesz resolvent operator for the
quaternionic functional calculus.
\begin{Dn}[The $S$-resolvent operator]
Let $V$ be a two sided quaternionic Banach space, $T\in\mathcal{B}(V)$ and $s\in
\rho_S(T)$.
We define the left $S$-resolvent operator as
\begin{equation}\label{quatSresolrddlft}
S_L^{-1}(s,T):=-(T^2-2 {\rm Re}\,(s) T+|s|^2\mathcal{I})^{-1}
(T-\overline{s}\mathcal{I}),
\end{equation}
and the right $S$-resolvent operator as
\begin{equation}\label{quatSresorig}
S_R^{-1}(s,T):=-(T-\overline{s}\mathcal{I})(T^2-2 {\rm Re}\,(s)
T+|s|^2\mathcal{I})^{-1}.
\end{equation}
\end{Dn}

\begin{Tm}
Let $T\in\mathcal{B}(V)$ and let $s \in \rho_S(T)$. Then, the left $S$-resolvent
operator satisfies the equation
\begin{equation}\label{quatSresolrddlftequ}
S_L^{-1}(s,T)s-TS_L^{-1}(s,T)=\mathcal{I},
\end{equation}
while the right $S$-resolvent
operator satisfies the equation
\begin{equation}\label{quatSresorigequa}
sS_R^{-1}(s,T)-S_R^{-1}(s,T)T=\mathcal{I}.
\end{equation}
\end{Tm}
\begin{Dn}
Let $V$ be a two sided quaternionic Banach space,   $T\in\mathcal{B}(V)$
and let $U \subset \mathbb{H}$ be an axially symmetric s-domain
that contains  the $S$-spectrum $\sigma_S(T)$ and such that
$\partial (U\cap \mathbb{C}_I)$ is union of a finite number of
continuously differentiable Jordan curves  for every $I\in\mathbb{S}$.
We say that $U$ is a $T$-admissible open set.
\end{Dn}
We can now introduce the class of functions for which we can define the two
versions of the quaternionic functional calculus.
 \begin{Dn}\label{quatdef3.9}
Let $V$ be a two sided quaternionic Banach space,   $T\in\mathcal{B}(V)$
and let  $W$ be an open set in $\hh$.
\begin{itemize}
\item[(1)]
A function  $f\in \mathcal{R}^L(W)$  is said to be locally left
hyperholomorphic  on $\sigma_S(T)$
if there exists a $T$-admissible domain $U\subset \hh$ such that
$\overline{U}\subset W$, on
which $f$ is left hyperholomorphic.
We will denote by $\mathcal{R}^L_{\sigma_S(T)}$ the set of locally
\index{$\mathcal{R}^L_{\sigma_S(T)}$}
left hyperholomorphic functions on $\sigma_S (T)$.
\item[(2)]
A function $f\in \mathcal{R}^R(W)$ is said to be locally right
hyperholomorphic on $\sigma_S(T)$
if there exists a $T$-admissible domain $U\subset \hh$ such that
$\overline{U}\subset W$, on
which $f$ is right hyperholomorphic.
We will denote by $\mathcal{R}^R_{\sigma_S(T)}$ the set of locally
\index{$\mathcal{R}^R_{\sigma_S(T)}$}
right hyperholomorphic functions on $\sigma_S (T)$.
\end{itemize}
\end{Dn}

Using the left $S$-resolvent operator $S_L^{-1} $, we now give
a result that motivates the functional
calculus; analogous considerations can be done using $S_R^{-1}$
with obvious modifications.
\begin{Dn}[The quaternionic functional calculus]\label{quatfunccalleftright}
Let $V$ be a two sided quaternionic Banach space and  $T\in\mathcal{B}(V)$.
  Let $U\subset \hh$ be a $T$-admissible domain and set $ds_I=- ds I$. We define
\begin{equation}\label{quatinteg311def}
f(T)={{1}\over{2\pi }} \int_{\partial (U\cap \mathbb{C}_I)} S_L^{-1} (s,T)\
ds_I \ f(s), \ \ {\it for}\ \ f\in \mathcal{R}^L_{\sigma_S(T)},
\end{equation}
and
\begin{equation}\label{quatinteg311rightdef}
f(T)={{1}\over{2\pi }} \int_{\partial (U\cap \mathbb{C}_I)} \  f(s)\ ds_I \
S_R^{-1} (s,T),\ \  {\it for}\ \ f\in \mathcal{R}^R_{\sigma_S(T)}.
\end{equation}
\end{Dn}

\section{Projectors, right and S-spectrum}
\setcounter{equation}{0}
\label{4}
An important result that we will prove in this section is
that the Riesz projector associated to a given quaternionic operator $T$ commute with $T$ itself.
We begin by recalling the definition of projectors and
some of their basic properties that still hold in the quaternionic setting.
\begin{Dn}
Let $V$ be a quaternionic Banach space. We say that $P$ is a projector if $P^2=P$.
\end{Dn}
It is easy to show that the following properties hold:
\begin{enumerate}
\item[(1)]
The range of $P$, denoted by $\ran(P)$ is closed.
\item[(2)]
$v\in \ran(P)$ if and only if $Pv=v$.
\item[(3)]
If $P$ is a projector also $I-P$ is a projector and $\ran(I-P)$ is closed.
\item[(4)]
$v\in \ran(I-P)$ if and only if $(I-P)v=v$, that is if and only
if $Pv=0$, as a consequence $\ran(I-P)=\ker(P)$.
\item[(5)]
For every $v\in V$ we have $v=Pv+(I-P)v$; $Pv\in \ran(P)$, $(I-P)v\in \ker(P)$.
So $v$ can be written as $v'=Pv$ and $v''=(I-P)v$. Since $\ran(P)\cap \ker(P)=\{0\}$
we have the decomposition
$V=\ran(P)\oplus \ker(P)$.
\end{enumerate}

\begin{Tm}\label{PTcommutation}
Let $T\in\mathcal{B}(V)$ and
 let $\sigma_S(T)= \sigma_{1S}\cup \sigma_{2S}$,
with ${\rm dist}\,( \sigma_{1S},\sigma_{2S})>0$. Let $U_1$ and
$U_2$ be two open sets such that  $\sigma_{1S} \subset U_1$ and $
\sigma_{2S}\subset U_2$,  with $\overline{U}_1
\cap\overline{U}_2=\emptyset$. Set
\begin{equation}\label{pigei}
P_j:=\frac{1}{2\pi }\int_{\partial (U_j\cap \mathbb{C}_I)}S_L^{-1}(s,T) \,
ds_I, \ \ \ \ \ j=1,2,
\end{equation}
\begin{equation}\label{tigei}
T_j:=\frac{1}{2\pi }\int_{\partial (U_j\cap \mathbb{C}_I)}S_L^{-1}(s,T) \,
ds_I\,s, \ \ \ \  j=1,2.
\end{equation}
Then the following properties hold:
\begin{itemize}
\item[(1)]
$P_j$ are projectors and $TP_j=P_jT$ for $j=1,2$.
\item[(2)] For $\lambda\in \rho_S(T)$ we have
\begin{equation}\label{tipigeieqL}
P_jS_L^{-1} (\lambda,T)\lambda-T_jS_L^{-1} (\lambda,T)=P_j, \ \ \ \ \ j=1,2,
\end{equation}
\begin{equation}\label{tipigeieqR}
\lambda S_R^{-1} (\lambda,T)P_j-S_R^{-1} (\lambda,T)T_j=P_j, \ \ \ \ \ j=1,2.
\end{equation}
\end{itemize}
\end{Tm}
\begin{proof}
The fact that $P_j$ are projectors is proved in \cite{MR2752913}.
Let us prove that $TP_j=P_jT$. Observe that the functions
$f(s)=s^m$, for $m\in \mathbb{N}_0$ are both right and left
 slice hyperholomorphic. So
 the operator $T$ can be written as
$$
T={{1}\over{2\pi }} \int_{\partial (U\cap \mathbb{C}_I)} S_L^{-1} (s,T)\  ds_I \ s=
{{1}\over{2\pi }} \int_{\partial (U\cap \mathbb{C}_I)} \  s\ ds_I \ S_R^{-1} (s,T);
$$
analogously, for the projectors $P_j$ we have
$$
P_j={{1}\over{2\pi }} \int_{\partial (U_j\cap \mathbb{C}_I)} S_L^{-1} (s,T)\
ds_I \ =
{{1}\over{2\pi }} \int_{\partial (U_j\cap \mathbb{C}_I)} \  \ ds_I \ S_R^{-1} (s,T).
$$
From the identity
$$
T_j={{1}\over{2\pi }} \int_{\partial (U_j\cap \mathbb{C}_I)} S_L^{-1} (s,T)\
ds_I \ s={{1}\over{2\pi }} \int_{\partial (U_j\cap \mathbb{C}_I)} \  s\ ds_I \
S_R^{-1} (s,T)
$$
we can compute $TP_j$ as:
$$
TP_j={{1}\over{2\pi }} \int_{\partial (U_j\cap \mathbb{C}_I)} TS_L^{-1} (s,T)\  ds_I \
$$
and using the resolvent equation (\ref{quatSresolrddlftequ})
it follows
$$
TP_j={{1}\over{2\pi }} \int_{\partial (U_j\cap \mathbb{C}_I)} [S_L^{-1}(s,T)\
s-\mathcal{I}]\  ds_I \ =
{{1}\over{2\pi }} \int_{\partial (U_j\cap \mathbb{C}_I)} S_L^{-1}(s,T)\ s\  ds_I
$$
$$
=
{{1}\over{2\pi }} \int_{\partial (U_j\cap \mathbb{C}_I)} S_L^{-1}(s,T)\  ds_I\ s=T_j.
$$
Now consider
$$
P_jT={{1}\over{2\pi }} \int_{\partial (U_j\cap \mathbb{C}_I)}
\  \ ds_I \ S_R^{-1} (s,T)T
$$
and using the resolvent equation (\ref{quatSresorigequa})
we obtain
$$
P_jT={{1}\over{2\pi }} \int_{\partial (U_j\cap \mathbb{C}_I)} \
 \ ds_I \ [s \ S_R^{-1}(s,T)-\mathcal{I}]
={{1}\over{2\pi }} \int_{\partial (U_j\cap \mathbb{C}_I)} \
\ ds_I \ s\  S_R^{-1}(s,T)=T_j
$$
so we have the equality $P_jT=TP_j$.
To prove (\ref{tipigeieqL}), for $\lambda\in \rho_S(T)$, consider and compute
$$
 P_jS_L^{-1} (\lambda,T)\lambda={{1}\over{2\pi }}
 \int_{\partial (U_j\cap \mathbb{C}_I)} \  \ ds_I \
 S_R^{-1} (s,T)S_L^{-1} (\lambda,T)\lambda.
$$
Using the S-resolvent equation (\ref{quatSresolrddlftequ}) it
follows that
$$
 P_jS_L^{-1} (\lambda,T)\lambda={{1}\over{2\pi }} \int_{\partial (U_j\cap
 \mathbb{C}_I)} \  \ ds_I \ S_R^{-1} (s,T)[TS_L^{-1} (\lambda,T)+I]
$$
$$
={{1}\over{2\pi }} \int_{\partial (U_j\cap \mathbb{C}_I)} \  \ ds_I \ [S_R^{-1}
(s,T)T]S_L^{-1} (\lambda,T)+P_j.
$$
By the S-resolvent equation (\ref{quatSresorigequa}) we get
$$
 P_jS_L^{-1} (\lambda,T)\lambda={{1}\over{2\pi }} \int_{\partial
 (U_j\cap \mathbb{C}_I)} \ s \ ds_I \ S_R^{-1} (s,T)S_L^{-1} (\lambda,T)+P_j
$$
$$
 =T_jS_L^{-1} (\lambda,T)+P_j
$$
which is (\ref{tipigeieqL}).
Relation (\ref{tipigeieqR}) can be proved in an analogous way.

\end{proof}

In  analogy with the classical case, we will call the operator
$P_j$  Riesz projector.\\

Our next result, of independent interest, is the validity of the
decomposition of the $S$-spectrum which is based on the Riesz
projectors. A simple but crucial result will be the following
Lemma:
\begin{La}\label{lieiruyh}
Let $T\in\mathcal{B}(V)$ and let $\lambda\in \rho_S(T)$.
Then the operator $(T^2-2\lambda_0T+|\lambda|^2\mathcal{I})^{-1}$
commutes with every operator $A$ that commutes with $T$.
\end{La}
\begin{proof}
Since $A$ commutes with $T$ we have that
$$
(T^2-2\lambda_0T+|\lambda|^2\mathcal{I})A=A(T^2-2\lambda_0T+|\lambda|^2\mathcal{I}).
$$
We get the statement by multiplying  the above relation on both sides by
$(T^2-2\lambda_0T+|\lambda|^2\mathcal{I})^{-1}$.
\end{proof}
Note that, unlike what happens in the classical case in which an operator $A$ which commutes with $T$  also commute
with the resolvent operator, here an operator $A$ commuting with $T$ just commute
with $(T^2-2\lambda_0T+|\lambda|^2\mathcal{I})^{-1}$. But this result is enough to prove the validity of the next theorem.

\begin{Tm} Let $T\in\mathcal{B}(V)$, suppose that $P_1$ is a
projector in $\mathcal{B}(V)$ commuting with $T$ and let $P_2=I-P_1$.
 Let  $V_j=P_j(V)$, $j=1,2$
and define the operators $T_j=TP_j=P_jT$. Denote by
$\widetilde{T}_j$ the restriction of $T_j$ to $V_j$, $j=1,2$.
Then
$$
\sigma_S(T)=\sigma_S(\widetilde{T}_1)\cup\sigma_S(\widetilde{T}_2).
$$
\end{Tm}
\begin{proof}
First of all note that $T=T_1+T_2$, \ $T_1(V_2)=T_2(V_1)=\{0\}$ and
that $T_j(V_j)\subseteq V_j$.

We have to show that $\rho_S(T)=\rho_S(\widetilde{T}_1)\cap\rho_S(\widetilde{T}_2)$.
Let us assume that $\lambda\in \rho_S(T)$ and  consider the identity
\begin{equation}\label{RPP}
T^2-2\lambda_0T+|\lambda|^2\mathcal{I}=
(T^2-2\lambda_0T+|\lambda|^2\mathcal{I})(P_1+P_2)
\end{equation}
$$
=(T^2_1-2\lambda_0T_1+|\lambda|^2P_1)+(T^2_2-2\lambda_0T_2+|\lambda|^2P_2).
$$
If we set
$$
Q_\lambda(T):=(T^2-2\lambda_0T+|\lambda|^2\mathcal{I})^{-1}
$$
we have
\begin{equation}\label{RPP1}
Q_\lambda(T)=(P_1+P_2)Q_\lambda(T)(P_1+P_2)=P_1Q_\lambda(T)P_1+P_2Q_\lambda(T)P_2;
\end{equation}
in fact, by Lemma \ref{lieiruyh} and  by the relation $P_1P_2=P_2P_1=0$,
we deduce $$P_1Q_\lambda(T)P_2=P_2Q_\lambda(T)P_1=0.$$
We now multiply the identity (\ref{RPP}) by $Q_\lambda(T)$ on the left and
by (\ref{RPP1}) we obtain
$$
\mathcal{I}=(P_1Q_\lambda(T)P_1+P_2Q_\lambda(T)P_2)[(T^2_1-2
\lambda_0T_1+|\lambda|^2P_1)+(T^2_2-2\lambda_0T_2+|\lambda|^2P_2)].
$$
Using again Lemma \ref{lieiruyh} and  $P_1P_2=P_2P_1=0$ we obtain
\begin{equation}\label{ident}
\mathcal{I}=P_1Q_\lambda(T)P_1(T^2_1-2\lambda_0T_1+|\lambda|^2P_1)+
P_2Q_\lambda(T)P_2(T^2_2-2\lambda_0T_2+|\lambda|^2P_2).
\end{equation}
Let us set
$$
Q_{\lambda, j}(T):=P_jQ_\lambda(T)P_j, \ \ \ \ j=1,2.
$$
It is immediate to observe that
$$
Q_{\lambda, j}(T)(V_j)\subseteq V_j,\ \ \ \ j=1,2,
$$
and from (\ref{ident}) we deduce
$$
Q_{\lambda, j}(T)(T^2_j-2\lambda_0T_j+|\lambda|^2P_j)=P_j, \ \ \ \ j=1,2.
$$
As a consequence, $Q_{\lambda, j}(T)$ restricted to $V_j$ is the inverse of
$(\widetilde{T}^2_j-2\lambda_0\widetilde{T}_j+|\lambda|^2P_j)$ and so we conclude that
$\lambda\in \rho_S(\widetilde{T}_1)\cap \rho_S(\widetilde{T}_2)$.

Conversely, assume that $\lambda\in \rho_S(\widetilde{T}_1)\cap \rho_S(\widetilde{T}_2)$.
Let us set
$$
\widetilde{Q}_{\lambda, j}(T):=(\widetilde{T}^2_j-2\lambda_0\widetilde{T}_j+|\lambda|^2P_j)^{-1}
$$
and define
$$
\widetilde{Q}=P_1\widetilde{Q}_{\lambda, 1}(T)P_1+P_2\widetilde{Q}_{\lambda, 2}(T)P_2.
$$
We have
$$
\widetilde{Q}(T^2-2\lambda_0T+|\lambda|^2\mathcal{I})
=[P_1\widetilde{Q}_{\lambda, 1}(T)P_1+P_2\widetilde{Q}_{\lambda, 2}(T)
P_2](T^2-2\lambda_0T+|\lambda|^2\mathcal{I})
$$
$$
=P_1(\widetilde{T}^2_1-2\lambda_0\widetilde{T}_1+|\lambda|^2P_1)^{-1}
P_1(T^2-2\lambda_0T+|\lambda|^2\mathcal{I})
$$
$$
+
P_2(\widetilde{T}^2_2-2\lambda_0\widetilde{T}_2+|\lambda|^2P_2)^{-1}P_2
(T^2-2\lambda_0T+|\lambda|^2\mathcal{I})
$$
$$
=P_1+P_2=\mathcal{I}.
$$
Analogously $(T^2-2\lambda_0T+|\lambda|^2\mathcal{I})\widetilde{Q}=\mathcal{I}$.
So $\lambda \in \rho_S(T)$.
\end{proof}

In all our discussions on the functional calculus we have used the notion of $S$-spectrum.
However, in the literature, also
other types of spectra are used: the so-called left spectrum and the right spectrum.
In order to discuss the notion of right spectrum it is not necessary to assume that $V$ is a two sided linear space, so we will consider quaternionic right linear spaces.
We recall the following definition:

\begin{Dn}
Let $T:V\to V$ be a right linear quaternionic operator on a right
quaternionic  Banach space $V$.
We denote by $\sigma_R(T)$ the right spectrum  of $T$
that is
$
\sigma_R(T)=\{ s\in \mathbb{H}\  :\ \ Tv= vs \ {\it for\ } v\in V , \ v\not=0 \}.
$
\end{Dn}
As it has been widely discussed in the literature, one can also
define the left spectrum, i.e. the set of $s\in\mathbb{H}$ such that $Tv=sv$.
However, the notion of left spectrum is not very useful, see \cite{MR1333599}.
The  S-spectrum and the left spectrum are not, in general, related, see \cite{MR2752913}.
The right spectrum is more useful and more studied.
 It has a structure similar to the one of the S-spectrum,
 indeed whenever it contains an element $s$, it contains also the whole 2-sphere
$[s]$. However the operator $\mathcal{I} s-T$, where
$(\mathcal{I}s )(v):=vs$, is not a right linear operator; thus
the notion of right spectrum is not associated to a linear
resolvent operator and this represents a disadvantage since it
prevents to define a functional calculus. The following result,
see \cite{CS_CRAS}, states that the right spectrum coincides with
the S-spectrum and thus $\sigma_R(T)$ can now be related to the
linear operator $T^2-2s_0T+|s|^2\mathcal{I}$.
\begin{Tm}\label{S=R}
Let $T$ be a right linear quaternionic operator. Then its point $S$-spectrum coincides with
the right spectrum.
\end{Tm}
Theorem \ref{S=R} is crucial since all known results
on the right spectrum become valid also for the point S-spectrum.
\\
Let us now consider the two eigenvalue problems:
$$
Tv=v s, \ \ \ v\not=0,
$$
and
$$
(T^2-2s_0 T+|s|^2\mathcal{I})w=0, \ \ \ \ \ \ w\not=0.
$$
As it is well known the right eigenvectors do not form a right
linear subspace of $V$, while the
$S$-eigenvectors do, as it is immediate to verify.
We have the following proposition which will be useful in the sequel.
\begin{Pn}\label{eigenvector}
Let $v$ be a right eigenvector associated to $s\in \sigma_R(T)$. Then
we have
$$
(T^2-2 s_0 T+|s|^2\mathcal{I})v=0.
$$
\end{Pn}
\begin{proof}
Since $Tv=v s$ it follows that $T^2v=T(v s)=v s^2$.
Thus we have
$$
(T^2-2 \ s_0 T+|s|^2\mathcal{I})v=vs^2-2s_0
sv+|s|^2v=v(s^2-2s_0 s+|s|^2)=0
$$
where we have used the identity $s^2-2s_0 s+|s|^2=0$ which
holds for every $s\in \mathbb{H}$.
\end{proof}

\section{A result on negative squares}
\setcounter{equation}{0}
\label{6}
In this section we will consider power series of the form $K(p,q)=\sum_{n,m=0}^\infty p^na_{n,m}\overline{q}^m$,
where $a_{n,m}=a_{n,m}^*\in\mathbb{H}^{N\times N}$. It is immediate that
$K(p,q)$ is a function slice hyperholomorphic in $p$ and right slice hyperholomorphic in $\bar q$; moreover the assumption on the coefficients $a_{n,m}$ implies that $K(p,q)$ is hermitian.
\begin{Pn}
Let $(a_{n,m})_{n,m\in\mathbb N_0}$ denote a sequence of $N\times
N$ quaternionic matrices such that $a_{n,m}=a_{m,n}^*$, and
assume that the power series
\[
K(p,q)=\sum_{n,m=0}^\infty p^na_{n,m}\overline{q}^m
\]
converges in a neighborhood $V$ of the origin. Then the following
are
equivalent:\\
$(1)$ The function $K(p,q)$ has $\kappa$ negative squares in
$V$.\\
$(2)$ All the finite matrices $A_{\mu}\stackrel{\rm def.}{=}
(a_{n,m})_{n,m=0,\ldots \mu}$ have at most $\kappa$ strictly
negative eigenvalues, and exactly $\kappa$ strictly negative
eigenvalues for at least one $\mu\in\mathbb N_0$.
\label{pnneq}
\end{Pn}

\begin{proof}
Let $r>0$ be such that $B(0,r)\subset V$, and let $I,J$ be two
units in the unit sphere of purely imaginary quaternions $\mathbb
S$ (see \eqref{sphere} for the latter). Then
\[
a_{n,m}=\frac{1}{4r^{n+m}\pi^2}\iint_{[0,2\pi]^2}e^{-Int}K(re^{It},re^{Js})e^{Jms}dtds.
\]
This expression does not depend on the specific choice of $I$ and
$J$. Furthermore, we take $I=J$ and so:
\[
A_\mu=\frac{1}{4r^{n+m}\pi^2}\iint_{[0,2\pi]^2}
\begin{pmatrix}I_N\\e^{-Jt}I_N\\ \vdots \\ e^{-J\mu t}I_N\end{pmatrix}
K(re^{Jt},re^{Js})\begin{pmatrix}I_N&e^{Js}I_N& \cdots & e^{J\mu
s}I_N\end{pmatrix}dtds.
\]
Write now
\[
K(p,q)=K_+(p,q)-F(p)F(q)^*,
\]
where $F$ is $\mathbb H^{N\times \kappa}$-valued. The function
$F$ is built from functions of the form $p\mapsto K(p,q)$ for a
finite number of $q$'s, and so is a continuous function of $p$,
and so is $K_+(p,q)$. See \cite[pp. 8-9]{adrs}. Thus
\[
A_\mu=A_{\mu,+}-A_{\mu,-}
\]
where
\[
\begin{split}
A_{\mu, +}&= \frac{1}{4r^{n+m}\pi^2}\iint_{[0,2\pi]^2}
\begin{pmatrix}I_N\\e^{-Jt}I_N\\ \vdots \\ e^{-J\mu t}I_N\end{pmatrix}
K_+(re^{Jt},re^{Js})\begin{pmatrix}I_N&e^{Js}I_N& \cdots & e^{J\mu
s}I_N\end{pmatrix}dtds\\
A_{\mu,-}&= \frac{1}{4r^{n+m}\pi^2}\iint_{[0,2\pi]^2}
\begin{pmatrix}I_N\\e^{-Jt}I_N\\ \vdots \\ e^{-J\mu t}I_N\end{pmatrix}
F(re^{Jt})F(re^{Js})^*\begin{pmatrix}I_N&e^{Js}I_N& \cdots &
e^{J\mu s}I_N\end{pmatrix}dtds.
\end{split}
\]
These expressions show that $A_\mu$ has at most $\kappa$ strictly
negative eigenvalues.\\

Conversely, assume that all the matrices $A_\mu$ have at most
$\kappa$ strictly negative eigenvalues, and define
\[
K_\mu(p,q)=\sum_{n,m=0}^\mu p^ma_{n,m}\overline{q}^m.
\]
Then, $K_\mu$ has at most $\kappa$ negative squares, as is seen by
writing $A_\mu$ as a difference of two positive matrices, one of
rank $\kappa$. Since, pointwise,
\[
K(p,q)=\lim_{\mu\rightarrow\infty} K_\mu(p,q),
\]
the function $K(p,q)$ has at most $\kappa$ negative squares.\\

To conclude the proof, it remains to see that the number of
negative squares of $K(p,q)$ and $A_\mu$ is the same. Assume that
$K(p,q)$ has $\kappa$ negative squares, but that the $A_\mu$ have
at most $\kappa^\prime<\kappa$ strictly negative eigenvalues.
Then, the argument above shows that $K(p,q)$ would have at most
$\kappa^\prime$ negative squares, which contradicts the
hypothesis. The other direction is proved in the same way.
\end{proof}

As consequences we have:

\begin{Pn}
In the notation of the preceding proposition, the number of
negative squares is independent of the neighborhood $V$.
\end{Pn}
\begin{proof}
This is because the coefficients $a_{n,m}$ do not depend on the
given neighborhood.
\end{proof}

\begin{Pn}
\label{pn51} Assume that $K(p,q)$  is $\mathbb H^{N\times
N}$-valued and has $\kappa$ negative squares in $V$ and let
$\alpha(p)$ be a $\mathbb H^{N\times N}$-valued slice
hyperholomorphic function and such that $\alpha(0)$ is
invertible. Then the function
\begin{equation}
\label{aka}
B(p,q)=\alpha(p)\star K(p,q)\star_r\alpha(q)^*
\end{equation}
has $\kappa$ negative squares in $V$.
\end{Pn}

\begin{proof}
Write $K(p,q)=\sum_{n,m=0}^\infty p^na_{n,m}\overline{q}^m$ and
$\alpha(p)=\alpha_0+p\alpha_1+\cdots$. The $\mu\times \mu$ main
block matrix $B_\mu$ corresponding to the power series
\eqref{aka} is equal to
\[
B_\mu=LA_\mu L^*,
\]
where
\[
L=\begin{pmatrix} \alpha_0&0&0&\cdots &0\\
                   \alpha_1&\alpha_0&0&\cdots&0\\
                   \alpha_2&\alpha_1&\alpha_0&0&\cdots\\
                   \vdots&\vdots& & &\\
                   \alpha_\mu&\alpha_{\mu-1}&\cdots &\alpha_1&\alpha_0
                   \end{pmatrix}
\]
Since $\alpha_0=\alpha(0)$ is assumed invertible, the signatures of $A_\mu$
and $B_\mu$ are the same for every $\mu\in\mathbb N_0$. By
Proposition \ref{pnneq} it follows that the kernels $A$ and $B$
have the same number of negative squares.
\end{proof}

\section{Operators in quaternionic Pontryagin spaces}
\setcounter{equation}{0}
\label{7}
This section contains some definitions and results on right quaternionic
 Pontryagin spaces. Some of the statements hold when we replace  Pontryagin spaces by
Krein spaces. \\

The following result, proved in the complex plane in
\cite[Theorem 2.4, p. 18]{ikl}, is very useful to study
convergence of sequences in Pontryagin spaces. It implies in
particular that in a reproducing kernel Pontryagin space,
convergence is equivalent to convergence of the self-inner product
together with pointwise convergence. The proof of the
quaternionic case appears in \cite[Proposition 12.9, p. 471]{as3}.

\begin{Pn}
\label{pn:ikl}
Let $(\mathscr P,[\cdot,\cdot])$ denote a quaternionic right
Pontryagin space. The sequence $f_n$ of elements in $\mathscr P$
tends to $f\in\mathscr P$ if and only if the following two
conditions hold:
\[
\begin{split}
\lim_{n\rightarrow\infty} [f_n,f_n]&=[f,f], \intertext{and}
\lim_{n\rightarrow\infty}[f_n,g]&=[f,g]\quad\mbox{for $g$ in a
dense subspace of $\mathscr P$.}
\end{split}
\]
\end{Pn}

We endow $\mathbb H^N$ with the inner product
\[
[u,v]_{\mathbb H^N}=v^*u.
\]
Furthermore, a Hermitian form will be defined as having the
following linearity condition:
\begin{equation}
\label{eqherm}
[fa,gb]=\overline{b}[f,g]a.
\end{equation}
\begin{Rk}{\rm When we consider two sided Pontryagin vector spaces, we require an additional property on the inner product
with respect to the left multiplication, i.e.
$$
[av,av]=|a|^2[v,v].
$$
This property is satisfied, for example, in $\mathbb{H}^N$ with the inner product described above.}
\end{Rk}

\begin{Tm}
\label{tm:milano}
Let $T$ be a contraction in a two sided quaternionic
Pontryagin space such that
$T$ has no S-spectrum on the unit sphere and it satisfies
$$
[S_L^{-1}(\lambda,T)\lambda v, S_L^{-1}(\lambda,T)\lambda v] \leq
[S_L^{-1}(\lambda,T)v, S_L^{-1}(\lambda,T)v], \ \ \ for \ \ \ |\lambda|=1.
$$
Then $T$ has maximal negative invariant subspace and this
subspace is unique.
\end{Tm}
\begin{proof}
Let $|\lambda|=1$ so that the operator  $S_L^{-1}(\lambda,T)$ exists.
The fact that $T$ is a contraction implies the inequality
$$
[TS_L^{-1}(\lambda,T)v, TS_L^{-1}(\lambda,T)v] <
[S_L^{-1}(\lambda,T)v, S_L^{-1}(\lambda,T)v]
$$
for $v\not=0$.
Using the $S$- resolvent equation one deduces
$$
[S_L^{-1}(\lambda,T) \lambda v+\mathcal{I}v,
S_L^{-1}(\lambda,T) \lambda v+\mathcal{I}v] <
[S_L^{-1}(\lambda,T)v, S_L^{-1}(\lambda,T)v]
$$
from which one gets
$$
 [S_L^{-1}(\lambda,T) \lambda v, S_L^{-1}(\lambda,T) \lambda v]
 + [v, v]
 + [S_L^{-1}(\lambda,T) \lambda v, v]
 +
 [ v, S_L^{-1}(\lambda,T)\lambda v]
 $$
 $$
 < [S_L^{-1}(\lambda,T)v, S_L^{-1}(\lambda,T)v].
$$
so, using the hypothesis, we finally get
$$
[v, v]+ [S_L^{-1}(\lambda,T) \lambda v, v]+[ v, S_L^{-1}(\lambda,T)
\lambda v]< 0.
$$
In the above inequality we replace $S_L^{-1}
(\lambda,T) \lambda$ by $S_L^{-1}(\lambda,T) \lambda d\lambda_I$, where $d
\lambda_I=-Ie^{I\theta}d\theta$
and integrate over $[0,2\pi]$.
Recalling the definition of Riesz projector
$$
P=-\frac{1}{2\pi}\int_{\partial(\mathbb{B}\cap \mathbb{C}_I)}
S_L^{-1}(\lambda,T) \lambda d\lambda_I
$$
we obtain
$$
[v, v]<   [P v, v]+[ v, Pv]
$$
and so
$$
[v, v]<  2 {\rm Re}\,[P v, v].
$$
 Theorem \ref{PTcommutation} implies that $PT=TP$ and the rest of
 the proof follows as in Theorem 11.1 p.76 in
 \cite{ikl}.
\end{proof}
For right quaternionic Pontryagin spaces we have the following result.
\begin{Pn}\label{maximal}
A contraction $T$ in a right quaternionic
Pontryagin space $\mathcal P$ possessing an eigenvalue $\lambda$ with $|\lambda|>1$  has maximal negative invariant subspace.
\end{Pn}
\begin{proof}  Let $v\not= 0$ be an eigenvector associated to the right eigenvalue $\lambda$. Then we have
$$
[Tv,Tv]=[v\lambda,v\lambda]< [v,v],
$$
from which we deduce
$$
|\lambda|^2[v,v]< [v,v]
$$
and so $[v,v]<0$. Consider the right subspace $\mathcal M$ generated by $v$. Then any element in $\mathcal M$ is of the form $va$, $a\in\mathbb{H}$ and
$[va,va]<0$. The subspace $\mathcal M$ is invariant under the action of $T$, indeed $T(va)=T(v)a=v\lambda a$. Thus $\mathcal M$ is a negative invariant subspace
of $\mathcal P$.  Then $\mathcal M$ is maximal or it is contained in another negative invariant subspace $\mathcal M_1$ and iterating this procedure we obtain a chain of inclusions $\mathcal M\subset \mathcal M_1\subset\ldots $ which should end because $\mathcal P_-$ is finite dimensional.
\end{proof}

In view of Definition \ref{neqs} below, it is useful to recall the
following result (see \cite[Corollary 6.2, p. 41]{MR97h:15020}).

\begin{Pn}
An Hermitian matrix $H$ with entries in $\mathbb H$ is
diagonalizable, and its eigenvalues are real. Furthermore,
eigenvectors corresponding to different eigenvalues are
orthogonal in $\mathbb H^N$. Let $(t,r,s)$ denote the signature of
$H$. There exists an invertible matrix $U\in\mathbb H^{N\times
N}$ such that
\begin{equation}
H=U\begin{pmatrix}\sigma_{tr}&0\\
0&0_{s\times s}\end{pmatrix}U^*
\end{equation}
where $\sigma_{tr}=\begin{pmatrix}I_t&0\\0&-I_r\end{pmatrix}$.
\label{pn:hermite}
\end{Pn}

\begin{Dn}
\label{neqs}
Let $A$ be a continuous right linear operator from the
quaternionic Pontryagin space $\mathscr P$ into itself. We say
that $A$ has $\kappa$ negative squares if for every choice of
$N\in\mathbb N$ and of $f_1,\ldots, f_N\in\mathscr P$, the
Hermitian matrix $H\in\mathbb H^{N\times N}$ with $jk$ entry
equal to
\begin{equation}
\label{a1}
[Af_k,f_j]_{\mathscr P}
\end{equation}
has at most $\kappa$ strictly negative eigenvalues, and exactly
$\kappa$ strictly negative eigenvalues for some choice of
$N,f_1,\ldots, f_N$.
\end{Dn}

Note that the above definition is coherent with the right linearity
condition \eqref{eqherm}. If we replace the $f_k$ by
$f_kh_k$ where $h_k\in\mathbb H$, the new matrix has $jk$ entry
\[
[Af_kh_k,f_jh_j]_{\mathscr P}=\overline{h_j}[Af_k,f_j]_{\mathscr P}h_k,
\]
and so
\[
\left([Af_kh_k,f_jh_j]_{\mathscr P}\right)_{j,k=1,\ldots, N}=D^*
\left([Af_k,f_j]_{\mathscr P}\right)_{j,k=1,\ldots, N}D,
\]
with
\[
D={\rm diag}~(h_1,h_2,\ldots, h_N).
\]
In case of left linear operators, \eqref{eqherm} is then replaced by
\[
[af,bg]=b[f,g]\overline{a},
\]
and the roles of $j$ and $k$ have to be interchanged in \eqref{a1}.
This problem does not appear in the commutative case.\\

We point out the following notation. Let $T$ be a bounded linear
operator from the quaternionic right Pontryagin space $(\mathscr
P_1,[\cdot,\cdot]_{\mathscr P_1} )$ into the quaternionic right
Pontryagin space $(\mathscr P_2,[\cdot,\cdot]_{\mathscr P_2})$,
and let $\s_1$ and $\s_2$ denote two signature operators such that
$(\mathscr P_1,\langle\cdot,\cdot\rangle_{\mathscr P_1})$ and
$(\mathscr P_2,\langle\cdot, \cdot\rangle_{\mathscr P_2})$ are
right quaternionic Pontryagin spaces, where
\[
\langle x,y\rangle_{\mathscr P_j}=[ x,\s_jy]_{\mathscr P_j},\quad
j=1,2.
\]

 We denote by $T^{[*]}$ the adjoint of
$T$ with respect to the Pontryagin structure and by $T^*$ its
adjoint with respect to the Hilbert space structure. Thus,
\[
\begin{split}
[Tx,y]_{\mathscr P_2}&=\langle Tx,\s_2y\rangle_{\mathscr P_2}\\
&=\langle x,T^*\s_2y\rangle_{\mathscr P_1}\\
&=[x,\s_1T^*\s_2y]_{\mathscr P_1},
\end{split}
\]
and so, as is well known in the complex case,
\[
T^{[*]}=\s_1T^*\s_2\quad{\rm and} \quad T^*=\s_1T^{[*]}\s_2.
\]

We will denote $\nu_-(A)=\kappa$. When $\kappa=0$ the operator is
called positive.
\begin{Tm}
\label{tm:facto} Let $A$ be a bounded right linear self-adjoint operator from
the quaternionic Pontryagin space $\mathscr P$ into itself, which
has a finite number of negative squares. Then, there exists a
quaternionic Pontryagin space $\mathscr P_1$ with ${\rm
ind}_{\mathscr P_1}=\nu_-(A)$, and a bounded right linear operator $T$ from
$\mathscr P$ into $\mathscr P_1$ such that
\[
A=T^{[*]}T.
\]
\end{Tm}

\begin{proof} The proof follows that of \cite[Theorem 3.4, p.
456]{MR2576304}, slightly adapted to the present non commutative
setting. Since $A$ is Hermitian, the formula
\[
[Af,Ag]_A=[Af,g]_{\mathscr P}
\]
defines a Hermitian form on the range of $A$. Since
$\nu_-(A)=\kappa$, there exists $N\in\mathbb N$ and $f_1,\ldots,
f_N\in\mathscr P$ such that the Hermitian matrix $M$ with $\ell
j$ entry $[Af_j,f_\ell]_{\mathscr P}$ has exactly $\kappa$
strictly negative eigenvalues. Let $v_1,\ldots , v_\kappa$ be the
corresponding eigenvectors, with strictly negative eigenvalues
$\lambda_1,\ldots, \lambda_\kappa$. As recalled in Proposition
\ref{pn:hermite} $v_j$ and $v_k$ are orthogonal when
$\lambda_j\not=\lambda_k$. We can, and will, assume that vectors
corresponding to a given eigenvalue are orthogonal. Then,
\begin{equation}
v_s^*Mv_t=\lambda_t\delta_{ts},\quad t,s=1,\ldots, N.
\label{Fortho}
\end{equation}
In view of \eqref{eqherm}, and with
\[
v_t=\begin{pmatrix}v_{t1}\\ v_{t2}\\ \vdots\\ v_{tN}\end{pmatrix},\quad t=1,\ldots, N,
\]
we see that \eqref{Fortho} can be rewritten as
\[
[F_s,F_t]_A=\lambda_t\delta_{ts},\quad{\rm with}\quad
F_s=\sum_{k=1}^N Af_kv_{sk},\quad t,s=1,\ldots, N.
\]
The space $\mathscr M$ spanned by $F_1,\ldots, F_N$ is strictly
negative, and it has an orthocomplement in $({\rm Ran}~
A,[\cdot,\cdot]_A)$, say $\mathscr M^{[\perp]}$, which is a right
quaternionic pre-Hilbert space. The space ${\rm Ran}~A$ endowed
with the quadratic form
\[
\langle m+h,m+h\rangle_A=-[m,m]_A+[h,h]_{A},\quad m\in\mathscr
M,\,\, h\in\mathscr M^{[\perp]},
\]
is a pre-Hilbert space, and we denote by $\mathscr P_1$ its
completion. We note that $\mathscr P_1$ is defined only up to an
isomorphism of Hilbert space. We denote by $\iota$ the injection
from ${\rm Ran}~A$ into $\mathscr P_1$ such that
\[
\langle f,f\rangle_A=\langle \iota (f),\iota(f)\rangle_{\mathscr
P_1}.
\]
We consider the decomposition $\mathscr P_1=\mathscr
\iota(M)\oplus\iota(M)^{\perp}$, and endow $\mathscr P_1$ with the
indefinite inner product
\[
[\iota(m)+h,\iota(m)+h]_{\mathscr P_!}
=[m,m]_A+\langle h, h\rangle_{\mathscr P_1}.
\]
See \cite[Theorem 2.5, p. 20]{ikl} for the similar argument in
the complex case. Still following \cite{MR2576304} we define
\[
Tf=\iota (Af),\quad f\in\mathscr P.
\]
We now prove that that $T$ is a bounded right linear operator from $\mathscr
P$ into $\iota({\rm Ran}~A)\subset\mathscr P_1$. Indeed, let
$(f_n)_{n\in\mathbb N}$ denote a sequence of elements in $\mathscr
P$ converging (in the topology of $\mathscr P$) to $f\in\mathscr
P$. Since ${\rm Ran A}$ is dense in $\mathscr P_1$, using
Proposition \ref{pn:ikl}  it is therefore enough
to prove that:\\

\[
\begin{split}
\lim_{n\rightarrow}[Tf_n,Tf_n]_{\mathscr P_1}&=[Tf,Tf]_{\mathscr P_1},\\
\intertext{{\rm and}}
\lim_{n\rightarrow\infty}[Tf_n,Tg]_{\mathscr P_1}&=[Tf,Tg]_{\mathscr P_1},
\quad \forall g\in\mathscr P.
\end{split}
\]
By definition of the inner product, the first equality amounts to
\[
\lim_{n\rightarrow}[Af_n,f_n]_{\mathscr P}=[Af,f]_{\mathscr P},
\]
which is true since $A$ is continuous, and similarly for the
second claim. Therefore $T$ has an adjoint operator, which is
also continuous. The equalities (with $f,g\in\mathscr P$)
\[
\begin{split}
[f,T^{[*]}Tg]_{\mathscr P}&=[Tf,Tg]_{\mathscr P_1}\\
                   &=[Tf,\iota(Ag)]_{\mathscr P_1}\\
                   &=[\iota(Af),\iota(Ag)]_{\mathscr P_1}\\
                   &=[Af,Ag]_{A}\\
                   &=[f,Ag]_{\mathscr P}
\end{split}
\]
show that $T^{[*]}T=A$.
\end{proof}

We note the following. As is well known, the completion of a
pre-Hilbert space is unique up to an isomorphism of Hilbert
spaces, and the completion need not be in general a subspace of
the original pre-Hilbert space. Some identification is needed. In
\cite{ikl} (see \cite[2.4, p. 19]{ikl} and also in
\cite{MR2576304}) the operator $\iota$ is not used, and the space
$\mathscr P_1$ is written directly as a direct sum of $\mathscr
M$ and of the completion of the orthogonal of $\mathscr M$. This
amounts to identify the orthogonal of $\mathscr M$ as a being a
subspace of its abstract completion.

\section{The structure theorem}
\setcounter{equation}{0}
\label{8}
We first give some background to provide motivation for the results presented in this section.
Denote by $R_0$ the backward-shift operator:
\[
R_0f(z)=\frac{f(z)-f(0)}{z}
\]
Beurling's theorem can be seen as the characterization of $R_0$-invariant subspaces of the Hardy space $\mathbf H_2(\mathbb D)$,
where $\mathbb D$ is the unit disk in $\mathbb{C}$. These
are the spaces $\mathbf H_2(\mathbb D)\ominus j\mathbf H_2(\mathbb D)$, where $j$ is an inner function. Equivalently, these
are the reproducing kernel Hilbert spaces with reproducing kernel $k_j(z,w)=\frac{1-j(z)\overline{j(w)}}{1-z\overline{w}}$
with $j$ inner. When replacing $j$ inner by $s$ analytic and contractive in the open unit disk, it is more difficult to
characterize reproducing kernel Hilbert spaces $\mathscr H(s)$
with reproducing kernel $k_s(z,w)$. Allowing for $s$ not necessarily scalar valued,
de Branges gave a characterization of $\mathscr H(s)$ spaces in \cite[Theorem 11, p. 171]{db-fact}.
This result was
extended in  \cite[Theorem 3.1.2, p. 85]{adrs} to the case of Pontryagin spaces.
The theorem below is the analog of de Branges' result in the slice-hyperholomorphic setting,
in which the backward-shift operator $R_0$ is now defined as
$$
R_0 f(p)=p^{-1}(f(p)-f(0))=(f(p)-f(0))\star_\ell p^{-1}.
$$
In order to prove the result, we will be in need of a fact which is direct consequence of Lemma 3.6 in \cite{acs2}: if $f$, $g$ are  two left slice hyperholomorphic functions then
$$
(f\star_l g)^*=g^*\star_r f^*.
$$
\begin{Tm}
Let $\s\in\mathbb H^{N\times N}$ be a signature matrix, and let
$\mathscr M$ be a Pontryagin space of $\mathbb H^N$-valued
functions slice hyperholomorphic in a spherical neighborhood $V$ of the
origin, and invariant under the operator $R_0$. Assume moreover
that
\begin{equation}
\label{structure} [R_0f,R_0f]_{\mathscr M}\le[f,f]_{\mathscr
M}-f(0)^*\s f(0).
\end{equation}
Then, there exists a Pontryagin space $\mathscr P$ such that
${\rm ind}_-\mathscr P=\nu_-(\s)$ and a $\mathbf L(\mathscr P,
\mathbb H^N)$-valued slice hyperholomorphic function $S$ such
that the elements of $\mathscr M$ are the restrictions to $V$ of
the elements of $\mathscr P(S)$.
\end{Tm}

\begin{proof} We follow the proof in \cite[Theorem 3.1.2, p. 85]{adrs}.
Let ${\mathscr P_2}=\mathscr M\oplus \hh_\s$, and denote by $C$
the point evaluation at the origin.
We divide the proof into a number of steps.\\

STEP 1: {\sl Let $p\in V$ and $f\in \mathscr M$. Then,
\begin{equation}
\label{pointevalua} f(p)=C\star (I-pR_0)^{-\star}f.
\end{equation}
}

STEP 2: {\sl The reproducing kernel of $\mathscr M$ is given by
\[
K(p,q)=C\star (I-pR_0)^{-\star}\left(C\star (I-qR_0)^{-\star}\right)^*.
\]
}

STEP 3:  {\sl Let $E$ denote the operator
\[
E=\begin{pmatrix}R_0\\ C\end{pmatrix}:\quad \mathscr
M\longrightarrow\mathscr P_2.
\]
There exists a quaternionic
Pontryagin space $\mathscr P_1$ with ${\rm ind}_{\mathscr
P_1}=\nu_-(J)$, and a bounded right linear operator $T$ from $\mathscr M$ into
$\mathscr P_1$ such that
\begin{equation}
\label{factoE}
I_{\mathscr M}-EE^{[*]}=T^{[*]}T.
\end{equation}
}

Write (see \cite[(1.3.14), p. 26]{adrs})
\[
\begin{split}
\begin{pmatrix}I_{\mathscr M}&0\\ E&I_{\mathscr P_2}\end{pmatrix}
\begin{pmatrix} I_{\mathscr M}&0\\
0&I_{\mathscr P_2}-EE^{[*]}\end{pmatrix}\begin{pmatrix}I_{\mathscr M}&E^{[*]}\\
0&I_{\mathscr P_2}\end{pmatrix}&=\\
&\hspace{-2cm}=
\begin{pmatrix}
I_{\mathscr M}&E^{[*]}\\ 0&I_{\mathscr
P_2}\end{pmatrix}\begin{pmatrix} I_{\mathscr
M}-E^{[*]}E&0\\0&I_{\mathscr P_2}
\end{pmatrix}\begin{pmatrix}I_{\mathscr M}&0\\
E&I_{\mathscr P_2}\end{pmatrix}.
\end{split}
\]
Thus,
\begin{equation}
\label{ineq4} \nu_-(I_{\mathscr P_2}-EE^{[*]})+\nu_-(\mathscr
M)=\nu_- (I_{\mathscr M}-E^{[*]}E)+\nu_-(\mathscr P_2),
\end{equation}
and noting that $\nu_-(\mathscr P_2)=\nu_-(\mathscr M)+\nu_-(\s)$,
we have (see also \cite[Theorem 1.3.4(1), p. 25]{adrs})
\[
\nu_-(I_{\mathscr P_2}-EE^{[*]})+\nu_-(\mathscr M)=\nu_-
(I_{\mathscr M}-E^{[*]}E)+\nu_-(\mathscr M)+\nu_-(\s).
\]
Equation \eqref{structure} can be rewritten as $I-E^{[*]}E\ge 0$, and
in particular $\nu_-(I-E^{[*]}E)=0$. Thus
\[
\nu_-(I_{\mathscr P_2}-EE^{[*]})=\nu_-(\s).
\]
Applying Theorem \ref{tm:facto} we obtain the factorization \eqref{factoE}.\\

We set
\[
T^{[*]}=\begin{pmatrix}B\\ D\end{pmatrix}:\,\, \mathscr
P_1\longrightarrow\mathscr M\oplus\hh_\s,
\]
and
\[
V=\begin{pmatrix}R_0&B\\ C&D\end{pmatrix}.
\]
Let
\[
S(p)=D+p C\star(I_{\mathscr M}-p A)^{-\star} B.
\]

STEP 4: {\sl We have that
\[
\s-S(p)\s S(q)^*=C\star(I-pR_0)^{-\star}\star (I-p\overline{q})
\s_{\mathscr M}((I-qA)^{-\star})^*\star_rC^*,
\]
where $\s_{\mathscr M}$ is a fundamental symmetry for $\mathscr
M$.}\\

The computation is as in our previous paper \cite{acs2}.
\end{proof}

We note that a corollary of \eqref{ineq4} is:

\begin{Tm}
Let $T$ be a contraction between right quaternionic Pontryagin
spaces of same index. Then, its adjoint is a contraction.
\end{Tm}

\begin{proof}
Indeed, when $\nu_-(\mathscr M)=\nu_-(\mathscr P_2)$ we have
\[
\nu_-(I_{\mathscr P_2}-EE^{[*]})=\nu_- (I_{\mathscr M}-E^{[*]}E).
\]
\end{proof}

\section{Blaschke products}
\setcounter{equation}{0}
\label{5}
As is well known and easy to
check, a rational function $r$ is analytic in the open unit disk
and takes unitary values on the unit circle if and only if it is
a finite Blaschke product. If one allows poles inside the unit
disk, then $r$ is a quotient of finite Blaschke products. This is a very special
case of a result of Krein and Langer discussed in Section \ref{9} below. In
particular, such a function cannot have a pole (or a zero) on the
unit circle. The case of matrix-valued rational functions which
take unitary-values (with respect to a possibly indefinite inner
product space) plays an important role in the theory of linear
systems. When the metric is indefinite, poles can occur on the
unit circle.
See for instance \cite{pootapov,gvkdm,bgr-ot64,ag}.\\

Slice hyperholomorphic functions have zeros that are either
isolated points or isolated 2-spheres. If a slice
hyperholomorphic function $f$ has zeros at $Z=\{a_1,a_2,\ldots,
[c_1], [c_2], \ldots \}$ then its reciprocal $f^{-\star}$ has
poles at the set $\{[a_1],[a_2],\ldots,  [c_1], [c_2], \ldots
\}$, $a_i,c_j\in\mathbb{H}$. So the sphere associated to a zero
of $f$ is a pole of $f^{-\star}$. In other words, the poles are
always 2-spheres as one may clearly see from the definition of
$f^{-\star}=(f^s)^{-1}f^c$, see also \cite{MR2572530}.\\

We now recall the definitions of
Blaschke factors, see \cite{acs2}, and then
discuss the counterpart of rational unitary functions
in the present setting. For the Blasckhe factors it is necessary to give two
different definitions, according to the fact that the zero of a
Blaschke factor is a point, see Definition \ref{Def5.2}, or a
sphere, see Definition \ref{Def5.13}.
\begin{Dn}\label{Def5.2}
Let $a\in\mathbb{H}$, $|a|<1$. The function
\begin{equation}
\label{eqBlaschke} B_a(p)=(1-p\bar
a)^{-\star}\star(a-p)\frac{\bar a}{|a|}
\end{equation}
is called  Blaschke factor at $a$.
 \end{Dn}

 \begin{Rk}{\rm
Let $a\in\mathbb{H}$, $|a|<1$. Then, see Theorem 5.5 in
\cite{acs2}, the Blaschke factor $B_a(q)$ takes the unit ball
$\mathbb{B}$ to itself and the boundary of the unit ball to
itself. Moreover, it has a unique zero for $p=a$.}
\end{Rk}
The Blaschke factor having zeros at a sphere is defined as
follows:
 \begin{Dn}\label{Def5.13}
Let $a\in\mathbb{H}$, $|a|<1$. The function
\begin{equation}
\label{blas_sph} B_{[a]}(p)=(1-2{\rm
Re}(a)p+p^2|a|^2)^{-1}(|a|^2-2{\rm Re}(a)p+p^2)
\end{equation}
is called  Blaschke factor at the sphere $[a]$.
 \end{Dn}
\begin{Rk}{\rm The definition of $B_{[a]}(p)$ does not depend on
the choice of the point $a$ that identifies the 2-sphere. In fact
all the elements in the sphere $[a]$ have the same real part and module.
It is immediate that the Blaschke factor $B_{[a]}(p)$ vanishes on
the sphere $[a]$.}
 \end{Rk}
The  following result has been proven in \cite{acs2}, Theorem
5.16:
\begin{Tm}
A Blaschke product having zeros at the set
 $$
 Z=\{(a_1,\mu_1), (a_2,\mu_2), \ldots, ([c_1],\nu_1), ([c_2],\nu_2), \ldots \}
 $$
 where $a_j\in \mathbb{B}$, $a_j$ have respective
 multiplicities $\mu_j\geq 1$, $a_j\not=0$ for
 $j=1,2,\ldots $, $[a_i]\not=[a_j]$ if $i\not=j$,
 $c_i\in \mathbb{B}$, the spheres $[c_j]$ have respective multiplicities $\nu_j\geq 1$,
 $j=1,2,\ldots$, $[c_i]\not=[c_j]$ if $i\not=j$
and
$$
\sum_{i,j\geq 1} \Big(\mu_i (1-|a_i|)+ \nu_j(1-|c_j|)\Big)<\infty
$$
is given by
\[
\prod_{i\geq 1} (B_{[c_i]}(p))^{\nu_i}\prod_{j\geq 1}^\star
(B_{a'_j}(p))^{\star \mu_j},
\]
where $a'_1=a_1$ and $a'_j\in [a_j]$, for $j=2,3,\ldots$,  are suitably chosen
elements.
\end{Tm}

\begin{Rk}{\rm It is not difficult to compute the slice hyperholomorphic inverses of the Blaschke factors using Definition \ref{reciprocal}.
The slice hyperholomorphic reciprocal of $B_a(p)$ and $B_{[a]}(p)$ are, respectively:
$$
B_a(p)^{-\star}=\frac{a}{|a|}(a-p)^{-\star}\star(1-p\bar
a),
$$
$$
B_{[a]}(p)^{-\star}=(|a|^2-2{\rm Re}(a)p+p^2)^{-1}(1-2{\rm
Re}(a)p+p^2|a|^2).
$$
The reciprocal of a Blaschke product is constructed by taking the reciprocal of the factors, in the reverse order.
}
\end{Rk}
\begin{Rk}{\rm The zeroes of $B_{[a]}(p)$ are poles of $B_{[a]}(p)^{-\star}$ and viceversa.
The Blaschke factor $B_{a}(p)$ has a zero at $p=a$ and a pole at the 2-sphere $[1/\bar a]$, while $B_{a}(p)^{-\star}$
has a zero at $p=1/\bar a$ and a pole at the 2-sphere $[a]$.
}
\end{Rk}
\begin{Pn}
Let ${\s}\in\mathbb H^{N\times N}$ denote a signature matrix (that
is, ${\s}={\s}^*={\s}^{-1}$) and let $(C,A)\in\mathbb H^{N\times
M}\times \mathbb H^{M\times M}$ be such that
$\cap_{n=0}^\infty\ker CA^n=\left\{0\right\}$. Let $P$ be an
invertible and Hermitian solution of the Stein equation
\begin{equation}
\label{eq:stein}
P-A^*PA=C^*{\s}C.
\end{equation}
Then, there exist matrices $(B,D)\in\mathbb H^{M\times N}\times
\mathbb H^{\times N\times N}$ such that the function
\begin{equation}
\label{realbgk}
S(p)=D+p C\star(I_M-pA)^{-\star}B
\end{equation}
satisfies
\begin{equation}
{\s}-S(p){\s}S(q)^*=C\star
(I_M-pA)^{-\star}(P^{-1}-pP^{-1}\overline{q})
\star_r(I_M-A^*\overline{q})^{-\star_r} \star_r C^*.
\end{equation}
\label{pn:bla}
\end{Pn}

Before the proof we mention the following. The vectors
$f_1,f_2,\ldots$ in the quaternionic Pontryagin space $(\mathcal
P, [\cdot, \cdot]_{\mathcal P})$ are said to be orthonormal if
\[
[f_j,f_\ell]_{\mathcal P}=\begin{cases}0,\,\,\,\,\quad{\rm if}\quad j\not=\ell,\\
\pm 1,\quad {\rm if}\quad j=\ell.\end{cases}
\]
The set $f_1,f_2,\ldots$  is called an orthonormal basis if the closed linear span of the
$f_j$ is all of $\mathcal P$. In the proof we used
the fact that in a finite dimensional quaternionic Pontryagin space
an orthonormal family can be extended to an orthonormal
basis. This is true because any non-degenerate closed space in
a quaternionic Pontryagin space admits an orthogonal complement. See
\cite[Proposition 10.3, p. 464]{as3}.\\

\begin{proof}[Proof of Proposition \ref{pn:bla}]
Following our previous paper \cite{acs2} the statement is
equivalent to find matrices $(B,D)\in\mathbb H^{M\times N}\times
\mathbb H^{\times N\times N}$ such that:
\begin{equation}
\label{convention}
\begin{pmatrix}A&B\\ C&D\end{pmatrix}\begin{pmatrix} P^{-1}&0\\
0&{\s}\end{pmatrix}\begin{pmatrix}A&B\\ C&D\end{pmatrix}^*=
\begin{pmatrix} P^{-1}&0\\
0&{\s}\end{pmatrix},
\end{equation}
or equivalently,

\begin{equation}
\begin{pmatrix}A&B\\ C&D\end{pmatrix}^*\begin{pmatrix} P&0\\
0&{\s}\end{pmatrix}\begin{pmatrix}A&B\\ C&D\end{pmatrix}=
\begin{pmatrix} P&0\\
0&{\s}\end{pmatrix}.
\end{equation}

By Proposition \ref{pn:hermite} there exists a matrix
$V\in\mathbb H^{M\times M}$ and $t_1,s_1\in\mathbb N_0$ such that
\[
P=V\s_{t_1s_1}V^*.
\]
Equation \eqref{eq:stein} can be then rewritten as
\[
V^{-1}A^*V\s_{t_1,s_1}V^*AV^{-*}+V^{-1}C^*{\s}CV^{-*}=\s_{t_1,s_1},
\]
and expresses that the columns of the $\mathbb H^{(M+N)\times M}$ matrix
\[
\begin{pmatrix}V^{*}AV^{-*}\\ CV^{-*}\end{pmatrix}=
\begin{pmatrix}V^{*}&0\\0&I_N\end{pmatrix}
\begin{pmatrix}A\\ C\end{pmatrix}V^{-*}
\]
are orthogonal in $\mathbb H^{M+N}$, endowed with the inner product
\begin{equation}
\label{innerprod} [u,v]=u_1^*\s_{t_1,s_1}u_1+u_2^*{\s}u_2,\quad
u=\begin{pmatrix}u_1\\ u_2\end{pmatrix},
\end{equation}
the first $t_1$ columns having self-inner product equal to $1$ and
the next $s_1$ columns having self-inner product equal to $-1$. We
can complete these columns to form an orthonormal basis of
$\mathbb H^{(M+N)\times(M+N)}$ endowed with the inner product
\eqref{innerprod}, that is, we find a matrix $X\in\mathbb
H^{(M+N)\times N}$
\[
\begin{pmatrix}
\begin{pmatrix}V^*AV^{-*}\\ CV^{-*}\end{pmatrix}&X\end{pmatrix}
\in\mathbb H^{(M+N)\times(M+N)}
\]
unitary with respect to \eqref{innerprod}. From
\[
\begin{pmatrix}
\begin{pmatrix}V^*AV^{-*}\\ CV^{-*}\end{pmatrix}&X\end{pmatrix}^*
\begin{pmatrix}\s_{t_1s_1}&0\\ 0&{\s}\end{pmatrix}
\begin{pmatrix}\begin{pmatrix}V^*AV^{-*}\\ CV^{-*}\end{pmatrix}&X\end{pmatrix}=
\begin{pmatrix}\s_
{t_1s_1}&0\\ 0&{\s}\end{pmatrix},
\]
we obtain \eqref{convention} with
\begin{equation}
\begin{pmatrix}
B\\ D\end{pmatrix}=X\begin{pmatrix}V^*&0\\ 0&I_N\end{pmatrix}.
\end{equation}
\end{proof}

When the signature matrix ${\s}$ is taken to be equal to $I_N$ we
can get another more explicit formula for $S$.

\begin{Pn}
In the notation and hypothesis of the previous theorem, assume
${\s}=I_N$. Then, $(I_M-A)$ is invertible and the function
\begin{equation}
\label{Scentered1} S(p)=I_N-(1-p)
C\star(I_M-pA)^{-\star}P^{-1}(I_M-A)^{-*}C^*
\end{equation}
satisfies
\begin{equation}
I_N-S(p)S(q)^*=C\star
(I_M-pA)^{-\star}(P^{-1}-pP^{-1}\overline{q})\star_r(I_M-A^*\overline{q})^{-\star_r}
\star_r C^*.
\end{equation}
\label{pn:bla2}
\end{Pn}

Note that formula \eqref{Scentered1} is not a realization of the
form \eqref{realbgk}. It can be brought to the form
\eqref{realbgk} by writing:
\[
\begin{split}
S(p)&=S(0)+S(p)-S(0)\\
&=I_N-CP^{-1}(I_M-A)^{-*}C^*+\\
&\hspace{5mm}+p
C\star(I_M-pA)^{-\star}(I_M-A)P^{-1}(I_M-A)^{-*}C^*.
\end{split}
\]
\begin{proof}[Proof of Proposition \ref{pn:bla2}]
We write for $p,q$ where the various expressions make sense
\[
\begin{split}
S(p)I_N S(q)^*-I_N&= (I_N-(1-p)
C\star(I_M-pA)^{-\star}P^{-1}(I_M-A)^{-*}C^*)\times\\
&\hspace{5mm}\times(I_N-(1-q)
C\star(I_M-qA)^{-\star}P^{-1}(I_M-A)^{-*}C^*)^*-I_N\\
&=C\star(I_M-pA)^{-\star}\star\Delta\star_r(I_M-A^*\overline{q})^{-\star_r}\star_r
C^*,
\end{split}
\]
where
\[
\begin{split}
\Delta&=-(1-p)
P^{-1}\star(I_M-A)^{-*}(I_M-A^*\overline{q})^{-\star_r}-(I_M-pA)\star
(I-\overline{q}(I_M-A)^{-1}P^{-1}\\
&\hspace{5mm}+P^{-1}(I_M-A)^{-*}C^*JC\star(1-p)\star_r(1-\overline{q})\star_r
(I_M-A)^{-1}P^{-1}.
\end{split}
\]
Taking into account the Stein equation \eqref{eq:stein} with
$\sigma=I_M$ we have:
\[
\Delta=P^{-1}(I_M-A)^{-*}\star\Delta_1\star_r (I_M-A)^{-1}P^{-1},
\]
where, after some computations,
\[
\begin{split}
\Delta_1&= \left\{-(1-p)\star(I_M-\overline{q}A^*)\star_r
P(I_M-A)-\right.\\
&\hspace{5mm}\left.-(I_M-A)^*P\star
(I_M-pA)\star_r(1-\overline{q})+(1-p)\star(P-A^*PA)\star_r
(1-\overline{q})\right\}\\
&=(I_M-A^*)P(I_M-A).
\end{split}
\]
\end{proof}

We note that a space $\mathscr P(S)$ can be finite dimensional
without $S$ being square. For instance
\[
S(p)=\frac{1}{\sqrt{2}}\begin{pmatrix}1&b_a(p)\end{pmatrix}.
\]
On the other hand, finite dimensional ${\mathscr P}(S)$ spaces for
square $S$ correspond to the ${\s}$-unitary functions studied in
linear system theory. The factorization theory of these functions
(that is, the slice-hyperholomorphic counterpart of
\cite{ad3,ag,ag2}) will be considered in a future publication.
\section{Krein-Langer factorization}
\setcounter{equation}{0}
\label{9}
In the classical case, functions $S$ for which the kernel
\eqref{skl} has a finite number of negative squares have a
special structure: they can be written as the quotient of a Schur
function and of a finite Blaschke product. This is a result of
Krein and Langer. See for instance \cite{kl1}. In this section we
present some related results.

\begin{Pn}
Let $S$ be a $\mathbb H^{N\times M}$-valued slice
hyperholomorphic function in $\mathbb B$ and of the form
\begin{equation}
\label{formlk}
S(p)=B(p)^{-\star}\star S_0(p),\quad p\in \mathbb B
\end{equation}
where $B$ is a $\mathbb H^{N\times N}$-valued Blaschke product and
$S_0$ is a $\mathbb H^{N\times M}$-valued Schur multiplier. Then,
$S$ is a generalized Schur function.
\end{Pn}
\begin{proof}
We follow the argument in \cite[\S 6.2]{ad3}. We have for
$n\in\mathbb N_0$ and $p,q\in \mathbb{B}$
\[
\begin{split}
p^n(I_N-S(p)S(q)^*)\overline{q}^n&=p^n
B(p)^{-\star}\star\left( B(p)B(q)^*-S_0(p)S_0(q)^*\right)
\star_r(B(q)^*)^{-\star_r}\overline{q}^n
\\
&=p^n B(p)^{-\star}\star\left(B(p)B(q)^*-I_N+\right.
\\
&\hspace{5mm}\left.+I_N-S_0(p)S_0(q)^*\right)
\star_r(B(q)^*)^{-\star_r}\overline{q}^n.
\end{split}
\]
Thus
\begin{equation}
\label{diff}
K_S(p,q)=B(p)^{-\star}\star\left(K_{S_0}(p,q)-K_B(p,q) \right)
(B(q)^*)^{-\star_r},
\end{equation}
where $K_{S_0}$ and $K_B$ are defined as in \eqref{skl}. Using
Proposition \ref{pn51} with $\kappa=0$ we see that formula
\eqref{diff} expresses the kernel $K_S$ as a difference of two
positive definite kernels, one being finite dimensional. It
follows that $K_S$ has a finite number of negative squares in
$\mathbb B$.
\end{proof}

\begin{Tm}
Let $S$ be a $\mathbb H^{N\times M}$-valued slice
hyperholomorphic function in $\mathbb B$, and such that the
associated space $\mathcal P(S)$ is finite dimensional. Then, $S$
admits a representation of the form \eqref{formlk}.
\end{Tm}
\begin{proof} Since the coefficients spaces are quaternionic
Hilbert spaces, $R_0$ is a contraction in $\mathcal P(S)$. We
proceed along the lines of \cite[\S 4.2 p. 141]{adrs}
and divide the proof in a number of steps.\\

STEP 1: {\sl The operator $R_0$ has no eigenvalues of modulus
$1$.}\\

Indeed, let $f\in\mathcal P(S)$ and $\lambda\in\mathbb H$ be such
that $R_0f=f\lambda$. Assume $|\lambda|=1$. From
\begin{equation}
\label{ineq}
[R_0f,R_0f]_{\mathcal P(S)}\le[f,f]_{\mathcal
P(S)}-f(0)^*f(0),\quad f\in\mathcal P(S)
\end{equation}
we get
\[
[f\lambda,f\lambda]\le [f,f]-f(0)^*f(0)
\]
and so $f(0)=0$. Reiterating \eqref{ineq} with $R_0f$ instead of
$f$ we get $(R_0f)(0)=0$, and in a similar way, $(R^n_0f)(0)=0$
for $n=2,3,\ldots$. But the $(R^n_0f)(0)$ are the coefficients of
the power series of $f$, and so $f=0$.\\

STEP 2: {\sl Let $\kappa$ be the number of negative  squares of
$K_S$. Then, $R_0$ has a $\kappa$-dimensional negative
invariant subspace.}\\

We write in matrix form $A=R_0$ and $C$ the point evaluation at
the origin, and denote by $H$ the matrix corresponding to the
inner product in $\mathcal P(S)$. Thus
\[
A^*HA\le H.
\]
Without loss of generality we assume that $A$ is in Jordan form
(see \cite{wiegmann}, \cite{MR97h:15020}) we denote by $\mathcal
L_+$ (resp. $\mathcal L_-$) the linear span of the generalized
eigenvectors corresponding to eigenvalues in $\mathbb B$ (resp.
outside the closure of $\mathbb B$). Since there are no
eigenvalues on $\partial\mathbb B$, $\mathbb H^N$ (where $N={\rm
dim}~\mathcal P(S)$) is spanned by $\mathcal L_+$ and $\mathcal
L_-$. As in \cite[Theorem 4.6.1, p. 57]{glr}, one shows that
\[
{\rm dim}~\mathcal L_+\le i_+(H)\quad{\rm and}\quad {\rm
dim}~\mathcal L_-\le i_-(H),
\]
where $i_+(H)$ is the number of positive eigenvalues of $H$ (and
similarly for $i_-(H)$), and by dimension argument equality holds
there. Thus $\mathcal L_-$ is a $\kappa$-dimensional invariant
subspace of $A$.\\

Let $G$ denote the solution of the matrix equation
\[
G-A^*GA=C^*C.
\]

STEP 3: {\sl  Let $\mathcal M$ be the space corresponding to
$\mathcal L_-$ in $\mathcal P(S)$, endowed with the metric
defined by $G$. Then $\mathcal M$ is
contractively included in $\mathcal P(S)$.}\\

Let $M$ denote the Gram matrix of $\mathcal M$ in the $\mathcal
P(S)$ inner product. We show that $M\ge P$. Indeed, in view of
\eqref{ineq}, the matrix $M$ satisfies
\[
A^*MA\le M-C^*C.
\]
In view of \eqref{eq:stein}, the matrix $M-P$ satisfies $
A^*(M-P)A\le M-P$, or equivalently (since $A$ is invertible)
\[
M-P\le A^{-*}(M-P)A^{-1}
\]
and so, for every $n\in\mathbb N$,
\begin{equation}
\label{eqpm}
M-P\le (A^{-*})^n(M-P)A^{-n}.
\end{equation}
Since the S-spectrum of $A$ is outside the closed unit ball, we
have by the S-spectral radius theorem (see \cite[Theorem 3.10, p. 616]
{MR2496568},\cite[Theorem 4.12.6, p. 155]{MR2752913}
\[
\lim_{n\rightarrow\infty}\|A^{-n}\|^{1/n}=0,
\]
and so $\lim_{n\rightarrow\infty}\|(A^{-*})^n(P-M)A^{-n}\|=0$.
Thus entrywise
\[
\lim_{n\rightarrow\infty}(A^{-*})^n(P-M)A^{-n}=0
\]
and it follows from \eqref{eqpm} that $M-P\le 0$.\\

By Proposition \ref{pn:bla},
\[
\mathcal M=\mathcal P(B).
\]
when $\mathcal M$ is endowed with the $P$ metric. Furthermore:\\

STEP 4: {\sl The kernel $K_S(p,q)-K_B(p,q)$ is positive.}\\

Let $k_{\mathcal M}(p,q)$ denote the reproducing kernel of
$\mathcal M$ when endowed with the $\mathcal P(S)$ metric. Then
\[
k_{\mathcal M}(p,q)-K_B(p,q)\ge 0
\]
and
\[
K_S(p,q)-k_{\mathcal M}(p,q)\ge 0.
\]
On the other hand
\[
K_S(p,q)-K_B(p,q)=K_S(p,q)-k_{\mathcal M}(p,q)+k_{\mathcal
M}(p,q)-K_B(p,q)
\]
and so is positive definite.\\

To conclude we apply Proposition \ref{pnneq} to
\[
K_S(p,q)-K_B(p,q)=B(p)\star\left(I_N-S_0(p)S_0(q)^*\right)\star_rB(q)^*
\]
where $S(p)=B(p)^{-\star}\star S_0(p)$, to check that $S_0$ is a
Schur function.
\end{proof}

\end{document}